# Global Stability of McKeithan's Kinetic Proofreading Model for T-Cell Receptor Signal Transduction

Eduardo D. Sontag

Department of Mathematics, Rutgers University

New Brunswick, NJ 08903

sontag@control.rutgers.edu

http://www.math.rutgers.edu/~sontag

March 7, 2018

**Abstract**

This note presents a self-contained and streamlined exposition of chemical-network results due to Feinberg, Horn, and Jackson. As an application, it shows the global asymptotic stability to equilibria of McKeithan's kinetic proofreading model for T-cell receptor signal transduction.

## 1   Introduction

This work was motivated a question which arose during Carla Wofsy's series of talks [4]. Consider the following system of first-order ordinary differential equations, for nonnegative functions $C_i(t)$:

$$
\begin{aligned}
\dot{C}_0 &= k_1 \left( T^* - \sum_{i=0}^{N} C_i \right) \left( M^* - \sum_{i=0}^{N} C_i \right) - (k_{-1,0} + k_{p,0}) C_0 \\
&\ \ \vdots \\
\dot{C}_i &= k_{p,i-1} C_{i-1} - (k_{-1,i} + k_{p,i}) C_i \\
&\ \ \vdots \\
\dot{C}_N &= k_{p,N-1} C_{N-1} - k_{-1,N} C_N
\end{aligned}
$$

(dots indicate derivatives with respect to time $t$) where the subscripted $k$'s, as well as $M^*$ and $T^*$, are arbitrary positive constants.

These equations represent the dynamics of the "kinetic proofreading" model proposed by McKeithan in [3] in order to describe how a chain of modifications of the T-cell receptor complex, via tyrosine phosphorylation and other reactions, may give rise to both increased sensitivity and selectivity of response. The quantities $C_i(t)$ represent concentrations of various intermediate complexes, and the assumption is that recognition signals are determined by the concentrations of the final complex $C_N$. The constant $k_1$ is the association rate constant for the reaction which produces an initial ligand-receptor complex $C_0$ from a T-cell receptor (TCR) and a peptide-major histocompatibility complex (MHC). The constants $k_{p,i}$ are the rate constants for each of the steps of phosphorylation or other intermediate modifications, and the constants $k_{-1,i}$ are dissociation rates. The differential equation for $C_0$ could also be written in an alternative manner, as

$$
\dot{C}_0 = k_1 TM - (k_{-1,0} + k_{p,0}) C_0,
$$

where $T(t)$ and $M(t)$ represent the concentrations of TCR's and MHC's respectively. The two conservation laws

$$
M + \sum_{i=0}^{N} C_i = M^*, \quad T + \sum_{i=0}^{N} C_i = T^*,
$$

when solved for $T$ and $M$, and substituted back into the equation, give the form shown earlier.

McKeithan's paper focused on the analysis of *equilibria* of the above differential equations (and made, for simplicity, the assumption that $k_{p,i} \equiv k_p$ and $k_{-1,i} \equiv k_{-1}$ for some fixed $k_p$ and $k_{-1}$). The question that arose during [4] was: what can be said about the *dynamics* of these equations?

We show here that the best possible conclusion is true: *there is a unique equilibrium point, and this equilibrium is globally asymptotically stable.*

It turns out that this result can be derived as an almost immediate corollary of the beautiful and powerful theory of *deficiency zero chemical reaction networks with mass-action kinetics* developed by Feinberg, Horn, and Jackson, cf. [1, 2]. In this note, we wish to provide a totally *self-contained and streamlined presentation* of the main theorems for such networks. In order to keep matters as simple as possible, however, we specialize to a subset ("single linkage class networks") which already contains the kinetic proofreading model; the reader is referred to [1], and the bibliography therein, for analogous results concerning more arbitrary deficiency zero networks.

The organization of this note is as follows. The results in question can be explained in terms of a special yet very general and appealing class of dynamical systems, which we introduce in Section 2, where the main theorems are also stated. In Section 3 we study equilibria and invariance properties for this class of systems. Section 4 has proofs of the stability theorems.

In Section 5, we specialize to the kinetic proofreading example. We do not actually use any terminology from chemical network theory, nor do we define the terms "deficiency zero" and "chemical reaction network," but also found in Section 5 are brief remarks concerning motivations for the general form of the systems considered, and the relation to chemical network concepts.

We wish to emphasize that this note is mainly expository, and credit for the mathematical results lies with Feinberg, Horn, and Jackson. It is our hope that this exposition will serve to make a wider audience in the dynamical systems and control theory communities aware of their work. A follow-up note [5] will add new results concerning robustness, dependence on parameters, and control-theoretic properties of the systems studied here.

**Acknowledgments**

The author wishes to express great thanks to Carla Wofsy for a fascinating series of lectures, to Leah Keshet for organizing a superb workshop, and most especially to Marty Feinberg for making available reprints of his work and for very enlightening e-mail discussions.

## 2  Definitions and Statements of Main Results

Some standard notations to be used are:

- $\mathbb{R}_{\geq 0}$ (resp., $\mathbb{R}_+$) = nonnegative (resp., positive) real numbers

- $\mathbb{R}_+^n$ (resp., $\mathbb{R}_+^{m \times m}$) = $n$-column vectors (resp., $m \times m$ matrices) with entries on $\mathbb{R}_+$; similarly for $\mathbb{R}_{\geq 0}$

- $\mathbb{R}_0^n$ = boundary of $\mathbb{R}_{\geq 0}^n$, set of vectors $x \in \mathbb{R}_{\geq 0}^n$ such that $x_i = 0$ for at least one $i \in \{1, \ldots, n\}$

- $x'$ = transpose of vector or matrix $x$

- $|x|$ = Euclidean norm of vector in $\mathbb{R}^n$

- $\langle x, z \rangle = x'z$, inner product of two vectors

- $\mathcal{D}^\perp = \{x \mid \langle x, z \rangle = 0 \; \forall z \in \mathcal{D}\}$.

The systems to be studied are parametrized by two matrices $A$ and $B$ with nonnegative entries, as well as a nonnegative function $\theta$, and have the following general form:

$$\dot{x} = f(x) = \sum_{i=1}^{m} \sum_{j=1}^{m} a_{ij}\, \theta(x_1)^{b_{1j}} \theta(x_2)^{b_{2j}} \ldots \theta(x_n)^{b_{nj}}\, (b_i - b_j) \tag{1}$$

where $b_\ell$ denotes the $\ell$-th column of $B$ (notice that the diagonal entries of $A$ are irrelevant, since $b_i - b_i = 0$). Several restrictions on $A$, $B$, and $\theta$ are imposed below. The powers are interpreted as follows, for any $r, c \geq 0$: $r^0 = 1$, $0^c = 0$ if $c > 0$, and $r^c = e^{c \ln r}$ if $r > 0$ and $c > 0$.

The main example of interest is when $\theta(y) = |y|$ and $B$ is a matrix whose entries are nonnegative integers. In that case, the equations (1) are polynomial for nonnegative vectors $x$. However, we will allow a more general class of systems.

We next describe the hypotheses on $\theta$, $A$, and $B$. The map

$$\theta\,:\, \mathbb{R} \to [0, \infty)$$

is locally Lipschitz, has $\theta(0) = 0$, satisfies $\int_0^1 |\ln \theta(y)|\, dy < \infty$, and its restriction to $\mathbb{R}_{\geq 0}$ is strictly increasing and onto. We suppose that

$$A = (a_{ij}) \in \mathbb{R}_{\geq 0}^{m \times m} \quad \text{is irreducible} \tag{2}$$

(that is, $(I + A)^{m-1} \in \mathbb{R}_+^{m \times m}$ or, equivalently, the incidence graph $G(A)$ is strongly connected, where $G(A)$ is the graph whose nodes are the integers $\{1, \ldots, m\}$ and for which there is an edge $j \to i$, $i \neq j$, if and only if $a_{ij} > 0$), and that

$$B = (b_1, \ldots, b_m) \in \mathbb{R}_{\geq 0}^{n \times m} \quad \text{has rank } m \tag{3}$$

(so, its columns $b_i$ are linearly independent),

$$\text{no row of } B \text{ vanishes,} \tag{4}$$

and

$$\text{each entry of } B \text{ is either } 0 \text{ or } \geq 1\,. \tag{5}$$

This last hypothesis insures that $f(x)$ in (1) is a locally Lipschitz vector field.

*From now on, we assume that all systems (1) considered satisfy the above assumptions.*

Our study will focus on those solutions of (1) which evolve in the nonnegative orthant $\mathbb{R}_{\geq 0}^n$. Recall that a subset $S \subseteq \mathbb{R}^n$ is said to be *forward invariant* with respect to the differential equation $\dot{x} = f(x)$ provided that each solution $x(\cdot)$ with $x(0) \in S$ has the property that $x(t) \in S$ for all positive $t$ in the domain of definition of $x(\cdot)$. We show in Section 3 that the nonnegative and positive orthants are forward invariant:

**Lemma 2.1** Both $\mathbb{R}_{\geq 0}^n$ and $\mathbb{R}_+^n$ are forward-invariant sets with respect to the system (1).

We will also show in Section 3 that there are no finite explosion times:

**Lemma 2.2** For each $\xi \in \mathbb{R}_{\geq 0}^n$ there is a (unique) solution $x(\cdot)$ of (1) with $x(0) = \xi$, defined for all $t \geq 0$.

In order to state concisely the main results for systems (1), we need to introduce a few additional objects. The subspace

$$\mathcal{D} := \operatorname{span} \{b_i - b_j,\, i \neq j\} = \operatorname{span} \{b_1 - b_2, \ldots, b_1 - b_m\} \tag{6}$$

can be seen as a distribution in the tangent space of $\mathbb{R}^n$; it has dimension $m - 1$ because adding $b_1$ to the last-shown generating set gives the column space of the rank-$m$ matrix $B$. For each vector $p \in \mathbb{R}^n$, we may also consider the parallel translate of $\mathcal{D}$ that passes through $p$, i.e. $p + \mathcal{D} = \{p + d, d \in \mathcal{D}\}$. A set $S$ which arises as an intersection of such an affine subspace with the nonnegative orthant:

$$S = (p + \mathcal{D}) \bigcap \mathbb{R}_{\geq 0}^n$$

(for some $p$, without loss of generality in $\mathbb{R}_{\geq 0}^n$) will be referred to as a *class*. If $S$ intersects the positive orthant $\mathbb{R}_+^n$, we say that $S$ is a *positive class*. The significance of classes is given by the fact that any solution $x(\cdot)$ of (1) must satisfy

$$x(t) - x(0) = \int_0^t f(x(s)) \, ds = \sum_{i=1}^m \sum_{j=1}^m \gamma(t) \, (b_i - b_j) \in \mathcal{D},$$

where $\gamma(t) = \int_0^t a_{ij} \theta(x_1(s))^{b_{1j}} \theta(x_2(s))^{b_{2j}} \ldots \theta(x_n(s))^{b_{nj}} \, ds$, so $x(t) \in x(0) + \mathcal{D}$ for all $t$. In particular:

**Lemma 2.3** Each class is forward invariant. □

We denote by $E$ (respectively, $E_+$ or $E_0$) the set of nonnegative (respectively, positive or boundary) equilibria of (1), i.e. the set of states $\bar{x} \in \mathbb{R}_{\geq 0}^n$ (respectively, $\in \mathbb{R}_+^n$ or $\in \mathbb{R}_0^n$), such that $f(\bar{x}) = 0$. Of course, $E$ is the disjoint union of $E_+$ and $E_0$.

**Theorem 1** *Consider any system (1), under the stated assumptions. For every maximal solution of (1) with $x(0) \in \mathbb{R}_{\geq 0}^n$, it holds that $x(t) \to E$ as $t \to +\infty$.*

This will be proved in Section 4, The invariance of classes (which are contained in subspaces of dimension $m - 1 < n$) precludes asymptotic stability of equilibria. The appropriate concept is that of *asymptotic stability relative to a class*. We say that an equilibrium $\bar{x} \in S$ is asymptotically stable relative to a class $S$ if it is (a) stable relative to $S$ (for each $\varepsilon > 0$, there is some $\delta > 0$ such that, for all solutions $x(\cdot)$, $|x(0) - \bar{x}| < \delta$ and $x(0) \in S$ imply $|x(t) - \bar{x}| < \varepsilon$ for all $t \geq 0$) and (b) locally attractive relative to $S$ (for some $\varepsilon > 0$, if $|x(0) - \bar{x}| < \varepsilon$ and $x(0) \in S$ then $x(t) \to \bar{x}$ as $t \to +\infty$). We say that $\bar{x} \in S$ is *globally* asymptotically stable relative to a class $S$ if it is stable relative to $S$ and globally attractive relative to $S$ ($x(t) \to \bar{x}$ for all solutions with $x(0) \in S$). The main results are as follows; the first part is shown in Section 3, and the remaining two in Section 4.

**Theorem 2** *Consider any system (1), under the stated assumptions. Fix any positive class $S$.*

a. *There is a unique equilibrium $\bar{x}_S \in S \bigcap E_+$.*

b. *The equilibrium $\bar{x}_S$ is asymptotically stable relative to $S$.*

c. *The equilibrium $\bar{x}_S$ is globally asymptotically stable relative to $S$ if and only if $S \bigcap E_0 = \emptyset$.*

**Example 2.4** The following trivial example may help in understanding the above theorems. We take $n = m = 2$, $\theta(y) = |y|$, $A =$ identity matrix, and

$$B = \begin{pmatrix} 1 & 2 \\ 1 & 1 \end{pmatrix}.$$

The system (1) is (for nonnegative states):

$$\dot{x}_1 = (x_1 - 1)x_1 x_2$$
$$\dot{x}_2 = 0$$

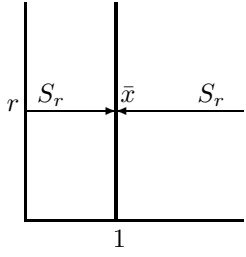

Figure 1: Example 2.4, dark lines indicate equilibria

and thus $E_0 = \mathbb{R}_0^2 = \{x \,|\, x_1 x_2 = 0\}$ and $E_+ = \{x \,|\, x_1 = 1, x_2 > 0\}$. The positive classes are the sets $S = S_r = \{x \,|\, x_1 \geq 0, x_2 = r\}$, for each $r > 0$, and for each such $S = S_r$, $\bar{x} = (1, r)'$ is asymptotically stable with domain of attraction $\{x \,|\, x_1 > 0, x_2 = r\}$. See Figure 1. Each class $S_r$ has a second equilibrium $(0, r)'$, but this second equilibrium is in the boundary, so there is no contradiction with part $a$ of Theorem 2. Regarding Theorem 1, observe that every trajectory either converges to an interior equilibrium $(1, r)'$ or it is itself a trajectory consisting of an equilibrium (and hence also converges to $E$, in a trivial sense). □

## 2.1 Other Expressions for the System Equations

The equations (1) have a considerable amount of structure, and various useful properties are reflected in alternative expressions for the system equations.

Let us introduce the map

$$\rho(y) := \ln \theta(y)$$

(with $\rho(0) = -\infty$). The restriction

$$\rho : \mathbb{R}_+ \to \mathbb{R}$$

is locally Lipschitz, strictly increasing, and onto $\mathbb{R}$. It also satisfies $\lim_{y \searrow 0} \rho(y) = -\infty$ and $\int_0^1 |\rho(y)| \, dy < \infty$. For any positive integer $n$, we let

$$\vec{\rho} : \mathbb{R}^n \to [-\infty, \infty)^n \, : \, x \mapsto (\rho(x_1), \dots, \rho(x_n))'$$

(we do not write "$\vec{\rho}_n$" to emphasize the dependence on $n$, because $n$ will be clear from the context). Then (1) can also be written as

$$\dot{x} = f(x) = \sum_{i=1}^m \sum_{j=1}^m a_{ij} e^{\langle b_j, \vec{\rho}(x) \rangle} \, (b_i - b_j) \, . \tag{7}$$

Here, the expression "$e^{\langle b_j, \vec{\rho}(x) \rangle}$" in (1) is interpreted in accordance with the conventions made for powers: if $x$ is a vector and $k \in \{1, \dots, n\}$ is an index such that $x_k = 0$ and $b_{kj} > 0$, then $e^{b_{kj}\rho(x_k)} = 0$, consistently with $e^{-\infty} = 0$, and thus also

$$e^{\langle b_j, \vec{\rho}(x) \rangle} = e^{b_{1j}\rho(x_1)} e^{b_{2j}\rho(x_2)} \cdots e^{b_{nj}\rho(x_n)} = 0 \, ,$$

but, if $b_{kj} = 0$, then we have $e^{b_{kj}\rho(x_k)} = 1$.

Another useful way of rewriting (1) is as follows. We write $f_k$ for the $k$-th coordinate of $f$ (i.e., the coordinates $x_k$ of solutions $x$ satisfy $\dot{x}_k = f_k(x)$). The terms in the sums defining $f_k$ can be collected into two disjoint sets: those that do not involve a product containing $\theta(x_k)$, for which $b_{kj} = 0$, and those which do involve $\theta(x_k)$. The latter, by assumption (5), have $b_{kj} \geq 1$, so we can factor $\theta(x_k)$ from $\theta(x_1)^{b_{1j}}\theta(x_2)^{b_{2j}}\ldots\theta(x_n)^{b_{nj}}$ and there remains a locally Lipschitz product. In other words, we can introduce, for each $k \in \{1,\ldots,n\}$, these two locally Lipschitz functions:

$$\alpha_k(x) := \sum_{j \in J_{k,1}} \left( \sum_{i=1}^m a_{ij} \ (b_{ki} - b_{kj}) \right) \theta(x_1)^{b_{1j}}\theta(x_2)^{b_{2j}}\ldots\theta(x_k)^{b_{kj}-1}\ldots\theta(x_n)^{b_{nj}} \qquad (8)$$

and

$$\beta_k(x) := \sum_{j \in J_{k,0}} \left( \sum_{i=1}^m a_{ij}b_{ki} \right) \theta(x_1)^{b_{1j}}\theta(x_2)^{b_{2j}}\ldots\theta(x_n)^{b_{nj}} , \qquad (9)$$

where $J_{k,1} := \{j \,|\, b_{kj} \geq 1\}$ and $J_{k,0} := \{j \,|\, b_{kj} = 0\}$. In terms of these,

$$f_k(x) = \alpha_k(x)\,\theta(x_k) + \beta_k(x) . \qquad (10)$$

In particular, since $\theta(0) = 0$,

$$x_k = 0 \ \Rightarrow \ f_k(x) = \beta_k(x) \qquad (11)$$

so, since $\beta_k(x) \geq 0$ for all $x$, the vector field $f$ always points towards the nonnegative orthant, on the boundary $\mathbb{R}_0^n$.

# 3  Equilibria and Invariance

We start with a global transversality result for $\mathcal{D}$ and $\mathcal{D}^\perp$. Then, we study interior and boundary equilibria, and provide basic results concerning the behavior of nonnegative solutions.

## 3.1  A Coordinatization Property

The next result shows provides, when specialized to classes, a one-to-one correspondence between points in $\mathbb{R}_+^n$ and pairs $(\mathcal{S}, \mathcal{R})$ consisting of a positive classes $\mathcal{S} = p + S$ and "coclasses" $\mathcal{R} = \vec{\rho}^{-1}(q + \vec{\rho}(S^\perp))$.

**Lemma 3.1** Let $D$ be any subspace of $\mathbb{R}^n$. For each $p, q$ in $\mathbb{R}_+^n$, there exists a unique $x = \varphi(p,q) \in \mathbb{R}_+^n$ such that:

$$x - p \in D \qquad (12)$$

and

$$\vec{\rho}(x) - \vec{\rho}(q) \in D^\perp . \qquad (13)$$

*Proof.* We start by introducing the following mapping, for each $i \in \{1,\ldots,n\}$:

$$L_i(t) := \int_0^{t+\rho(q_i)} \rho^{-1}(s)\,ds \, - \, p_i t$$

defined for $t \in \mathbb{R}$. Since $\rho^{-1}$ is an increasing onto map from $\mathbb{R}$ into $\mathbb{R}_+$, $L_i'(t) = \rho^{-1}(t + \rho(q_i)) - p_i$, and hence also $L_i(t)$, increases to infinity as $t \to +\infty$. Also, $L_i(t) \to +\infty$ as $t \to -\infty$, because $\rho^{-1}$ is nonnegative and $p_i > 0$. Thus, $L_i$ is proper, that is, $\{t \,|\, L_i(t) \leq v\}$ is compact for each $v$.

Now we take the (continuously differentiable) function

$$Q(y) := \sum_{i=1}^n L_i(y_i)$$

thought of as a function of $y \in \mathbb{R}^n$. This function is also proper, because

$$\{y \,|\, Q(y) \leq w\} \ \subseteq \ \prod_{i=1}^n \{t \,|\, L_i(t) \leq w - (n-1)\ell\}$$

where $\ell$ is any common lower bound for the functions $L_i$. Restricted to $D$, $Q$ is still proper, so it attains a minimum at some point $y \in D$. In particular, $y$ must be a critical point of $Q$ restricted to $D$, so

$$(\nabla Q(y))' = (\rho^{-1}(y_1 + \rho(q_1)) - p_1, \ldots, \rho^{-1}(y_n + \rho(q_n)) - p_n)' \in D^{\perp}. \tag{14}$$

Pick $x \in \mathbb{R}_+^n$ such that $\vec{\rho}(x) = y + \vec{\rho}(q)$. Then $\vec{\rho}(x) - \vec{\rho}(q) \in D$ by definition, and (14) gives also $x - p \in D^{\perp}$.

Finally, we show uniqueness. Suppose that there a second $z \in \mathbb{R}_+^n$ so that $z - p \in D$ and $\vec{\rho}(z) - \vec{\rho}(q) \in D^{\perp}$. This implies that $x - z \in D$ and $\vec{\rho}(x) - \vec{\rho}(z) \in D^{\perp}$. Since $\rho$ is an increasing function, we have that, for any two distinct numbers $a, b$, $(a - b)(\rho(a) - \rho(b)) > 0$. So

$$\sum_{i=1}^n (x_i - z_i)(\rho(x_i) - \rho(z_i)) = \langle x - z, \vec{\rho}(x) - \vec{\rho}(z) \rangle = 0.$$

implies $x = z$. ∎

The following quantity measures deviations relative to $\mathcal{D}^{\perp}$. Let us define, for each $x, z \in \mathbb{R}_+^n$:

$$\delta(x, z) := \sum_{i=1}^m \sum_{j=1}^m \left( \langle b_i - b_j, \vec{\rho}(x) - \vec{\rho}(z) \rangle \right)^2. \tag{15}$$

Note that $\delta(x, z) = 0$ if and only if $\vec{\rho}(x) - \vec{\rho}(z) \in \mathcal{D}^{\perp}$, since $\mathcal{D}$ is spanned by the differences $b_i - b_j$,

**Remark 3.2** We could also have defined a smaller, but basically equivalent, sum using only the generating differences $b_i - b_1$, $i = 1, \ldots, m - 1$, but the above definition for $\delta$ seems more natural. Moreover, note that if we let

$$\Delta(x, z) := \sum_{i=1}^{m-1} (\langle b_i - b_1, \vec{\rho}(x) - \vec{\rho}(z) \rangle)^2 + \sum_{\ell=m}^n (\langle v_\ell, x - z \rangle)^2,$$

where the $v_i$ constitute a basis of $\mathcal{D}^{\perp}$, then the uniqueness part of Lemma 3.1, applied with $D = \mathcal{D}$, gives that $\Delta(x, z) = 0$ if and only if $x = z$. (Because $x = \varphi(x, x)$, and $\Delta(x, z) = 0$ implies $z = \varphi(x, x)$.) □

## 3.2 Equilibria

It is convenient to also express the dynamics (1) in matrix terms. Letting

$$\widetilde{A} = A - \mathrm{diag}\left( \sum_{i=1}^m a_{i1}, \ldots, \sum_{i=1}^m a_{im} \right) = \begin{pmatrix} -\sum_{i \neq 1}^m a_{i1} & a_{12} & \cdots & a_{1m} \\ a_{21} & -\sum_{i \neq 2}^m a_{i2} & \cdots & a_{2m} \\ \vdots & \vdots & \vdots & \vdots \\ a_{m1} & a_{m2} & \cdots & -\sum_{i \neq m}^m a_{im} \end{pmatrix}$$

we write:

$$\dot{x} = f(x) = B\widetilde{A}\,\Theta_B(x) \tag{16}$$

where $\Theta_B$ is the mapping

$$\Theta_B : \mathbb{R}^n \to \mathbb{R}_{\geq 0}^m : x \mapsto \left( e^{\langle b_1, \vec{\rho}(x) \rangle}, \ldots, e^{\langle b_m, \vec{\rho}(x) \rangle} \right)'$$

obtained as the composition of the maps $x \mapsto \vec{\rho}(x)$, $z \mapsto B'z$, and $y \mapsto (e^{y_1}, \ldots, e^{y_m})'$. In particular, since rank $B = m$ and $\rho$ maps $\mathbb{R}_+$ onto $\mathbb{R}$,

$$\text{the restriction } \Theta_B : \mathbb{R}_+^n \to \mathbb{R}_+^m \text{ is onto.} \tag{17}$$

Note that $\Theta_B(x) = \left(\theta(x_1)^{b_{11}}\theta(x_2)^{b_{21}}\ldots\theta(x_n)^{b_{n1}},\ldots,\theta(x_1)^{b_{1m}}\theta(x_2)^{b_{2m}}\ldots\theta(x_n)^{b_{nm}}\right)'$.

For any two $\bar{x}, x \in \mathbb{R}^n_+$, it holds that:

$$(\exists\,\kappa > 0)\ \Theta_B(x) = \kappa\Theta_B(\bar{x}) \iff \vec{\rho}(x) - \vec{\rho}(\bar{x}) \in \mathcal{D}^\perp \tag{18}$$

(recall the definition (6) of $\mathcal{D}$). To see this, denote $y := \Theta_B(x)$, $\bar{y} := \Theta_B(\bar{x})$. If $y = \kappa\bar{y}$, then, with $k := \ln\kappa$, $\ln y_j = k + \ln\bar{y}_j$, for each $j = 1,\ldots,m$. Thus $\langle b_j, \vec{\rho}(x)\rangle = k + \langle b_j, \vec{\rho}(\bar{x})\rangle$, which implies that $\langle b_j, \vec{\rho}(x) - \vec{\rho}(\bar{x})\rangle = k$ for all $j$, and therefore

$$\langle b_j - b_i, \vec{\rho}(x) - \vec{\rho}(\bar{x})\rangle = 0 \quad \forall i, j\,.$$

Conversely, if this holds, we may define $k := \langle b_1, \vec{\rho}(x) - \vec{\rho}(\bar{x})\rangle$, $\kappa := e^k$, and reverse all implications.

We also note, using once again that $B$ has full column rank, that $f(\bar{x}) = 0$ is equivalent to $\widetilde{A}\,\Theta_B(\bar{x}) = 0$, that is:

**Lemma 3.3** A state $\bar{x}$ is an equilibrium if and only if $\Theta_B(\bar{x}) \in \ker\widetilde{A}$. $\qquad\qquad\square$

We now consider the matrix $\widetilde{A}$. The row vector $\underline{1} = (1,\ldots,1)$ has the property that $\underline{1}\widetilde{A} = (0,\ldots,0)$, so, in particular, $\widetilde{A}$ is singular. The following is a routine consequence of the Perron-Frobenius (or finite dimensional Krein-Rutman) Theorem.

**Lemma 3.4** There exists $\bar{y} \in \mathbb{R}^n_+ \bigcap \ker\widetilde{A}$ so that $\left(\mathbb{R}^n_{\geq 0} \setminus \{0\}\right)\bigcap\ker\widetilde{A} = \{\kappa\bar{y},\ \kappa > 0\}$.

*Proof.* If $y \in \mathbb{R}^n_{\geq 0}$ is any eigenvector of $\widetilde{A}$, corresponding to an eigenvalue $\lambda$, it follows that $0 = \underline{1}\widetilde{A}y = \underline{1}\lambda y = \lambda q$, where $q := \underline{1}y$ is a positive number (because $y$, being an eigenvector, is nonzero), and therefore necessarily $\lambda = 0$. In other words, a nonnegative eigenvector can only be associated to the zero eigenvalue. Pick now any $\gamma > 0$ large enough such that all entries of $\widehat{A} := \widetilde{A} + \gamma I$ are nonnegative. Since the incidence graph $G(\widehat{A})$ coincides with $G(A)$, it follows that $\widehat{A}$ is also irreducible. By the Perron-Frobenius Theorem, the spectral radius $\sigma$ of $\widehat{A}$ is positive and it is an eigenvalue of $\widehat{A}$ of algebraic multiplicity one, with an associated positive eigenvector $\bar{y} \in \mathbb{R}^n_+$. Moreover, every nonnegative eigenvector $y \in \mathbb{R}^n_{\geq 0}$ associated to $\sigma$ is a positive multiple of $\bar{y}$. As adding $\gamma I$ moves eigenvalues by $\gamma$ while preserving eigenvectors (that is, $(\widetilde{A} + \gamma I)y = (\lambda + \gamma)y$ is the same as $\widetilde{A}y = \lambda y$), $\bar{y}$ is a positive eigenvector of the original matrix $\widetilde{A}$. It is necessarily in the kernel of $\widetilde{A}$, since we already remarked that any nonnegative eigenvector must be associated to zero. Finally, if $y$ is any other nonnegative eigenvector of $\widehat{A}$, and in particular any element of $\left(\mathbb{R}^n_{\geq 0} \setminus \{0\}\right)\bigcap\ker\widetilde{A}$, then it is also a nonnegative eigenvector of $\widehat{A}$, and thus it must be a positive multiple of $\bar{y}$, completing the proof. $\qquad\blacksquare$

**Corollary 3.5** The set of positive equilibria $E_+$ is nonempty. Moreover, pick any fixed $\bar{x} \in E_+$. Then, for any positive vector $x \in \mathbb{R}^n_+$, the following equivalence holds:

$$x \in E_+ \iff \delta(x, z) = 0\,. \tag{19}$$

*Proof.* By Lemma 3.4, there is some $\bar{y} \in \mathbb{R}^n_+$ in $\ker\widetilde{A}$. By (17), there is some $\bar{x} \in \mathbb{R}^n_+$ such that $\Theta_B(\bar{x}) = \bar{y}$. In view of Lemma 3.3, $\bar{x}$ is an equilibrium.

Now fix any $\bar{x} \in E_+$ and any $x \in \mathbb{R}^n_+$, and let $y := \Theta_B(x)$, $\bar{y} := \Theta_B(\bar{x})$. Suppose $x \in E_+$. By Lemma 3.3, $y \in \ker\widetilde{A}$. By Lemma 3.4, every two positive eigenvectors of $\widetilde{A}$ are multiples of each other, so there is some $\kappa \in \mathbb{R}_+$ such that $y = \kappa\bar{y}$. By (18), $\vec{\rho}(x) - \vec{\rho}(\bar{x}) \in \mathcal{D}^\perp$. Conversely, if this holds, then, again by (18), $y = \kappa\bar{y}$. hence $y$ is also an eigenvector of $\widetilde{A}$, so by Lemma 3.3 we conclude $x \in E_+$. $\qquad\blacksquare$

### 3.3  Proof of Part a in Theorem 2

Pick any positive class $S = p + \mathcal{D}$, $p \in \mathbb{R}_+^n$. By Corollary 3.5, there is some equilibrium $\bar{x} \in E_+$. We apply Lemma 3.1 with $D = \mathcal{D}$, to obtain $\bar{x}_S = \varphi(p, \bar{x})$. By (13) and (19), $\bar{x}_S \in E_+$. By (12), $\bar{x}_S - p \in \mathcal{D}$, i.e., also $\bar{x}_S \in S$, as required. To show uniqueness, suppose that also $z \in S \bigcap E_+$. Since $z \in S$, $\bar{x}_S - z \in \mathcal{D}$, and since $z \in E_+$, $\vec{p}(\bar{x}_S) - \vec{p}(z) \in \mathcal{D}^\perp$ (by Corollary 19). Thus the uniqueness assertion in Lemma 3.1 gives $\bar{x}_S = z$. ∎

### 3.4  Boundary Equilibria

Fix any $x = (x_1, \ldots, x_n) \in \mathbb{R}_{\geq 0}^n$. We wish to study the implications of some coordinate $x_k$ vanishing. For each $j \in \{1, \ldots, m\}$ we consider the set

$$\mathcal{S}_j := \{ k \,|\, b_{kj} > 0 \}$$

which is nonempty, by (3). Note that $k \in \mathcal{S}_j$ if and only if $j \in J_{k,1}$.

We will use repeatedly the following fact, for any $j \in \{1, \ldots, m\}$ and $k \in \{1, \ldots, n\}$:

$$x_k = 0 \ \text{ and } \ k \in \mathcal{S}_j \ \Rightarrow \ \theta(x_1)^{b_{1j}} \theta(x_2)^{b_{2j}} \ldots \theta(x_n)^{b_{nj}} = 0 \tag{20}$$

which is obvious, since $\theta(x_k)^{b_{kj}}$ vanishes when $x_k = 0$ and $b_{kj} \neq 0$. In particular,

$$x_k = 0 \ \Rightarrow \ \theta(x_1)^{b_{1j}} \theta(x_2)^{b_{2j}} \ldots \theta(x_n)^{b_{nj}} b_{kj} = 0 \tag{21}$$

since either $b_{kj} = 0$ or (20) applies.

**Lemma 3.6** Take any $i, j \in \{1, \ldots, m\}$ such that $a_{ij} \neq 0$. Then,

$$(\forall \ell \in \mathcal{S}_j)\,(x_\ell > 0) \quad \Rightarrow \quad (\forall k \in \mathcal{S}_i)\,(x_k > 0 \text{ or } f_k(x) > 0)\,. \tag{22}$$

*Proof.* Pick any $k \in \mathcal{S}_i$, and assume that $x_k = 0$. The assumption "$\ell \in \mathcal{S}_j \Rightarrow x_\ell > 0$" is equivalent to $j \in J_{\ell,1} \Rightarrow x_\ell > 0$. So, since $x_k = 0$, necessarily $j \in J_{k,0}$. By (11), $f_k(x) = \beta_k(x)$, where $\beta$ is as in (9), and the index $j$ being considered does appear in the sum defining $\beta_k$. Moreover, for each $\ell \in \mathcal{S}_j$, $x_\ell \neq 0$ by hypothesis, so $\theta(x_1)^{b_{1j}} \theta(x_2)^{b_{2j}} \ldots \theta(x_n)^{b_{nj}} > 0$. On the other hand, $k \in \mathcal{S}_i$ means that $b_{ki} > 0$, and also $a_{ij} \neq 0$. Thus, the term involving this particular $i$ and $j$ in the sums defining $\beta_k(x)$ is positive (and the remaining terms are nonnegative). ∎

**Proposition 3.7** Take any $x \in \mathbb{R}_{\geq 0}^n$. The following properties are equivalent:

1. $x \in E_0$.

2. For every $j \in \{1, \ldots, m\}$ there is some $k \in \mathcal{S}_j$ such that $x_k = 0$.

3. For every $j \in \{1, \ldots, m\}$, $\theta(x_1)^{b_{1j}} \theta(x_2)^{b_{2j}} \ldots \theta(x_n)^{b_{nj}} = 0$.

*Proof.* [1 $\Rightarrow$ 2] Pick any $x \in E_0$. If the second property is false, then there is some index $j$ such that $x_\ell > 0$ for all $\ell \in \mathcal{S}_j$. We claim that for every index $i$, $x_k > 0$ for all $k \in \mathcal{S}_i$. Since $\bigcup_j \mathcal{S}_j = \{1, \ldots, n\}$ (recall hypothesis (4)), this will mean that $x_k > 0$ for all $k$, so $x$ could not have been a boundary point, a contradiction. Let $J = \{j \,|\, x_\ell > 0\, \forall \ell \in \mathcal{S}_j\}$ and let $I = \{1, \ldots, m\} \setminus J$. We know that $J \neq \emptyset$ and must prove that $I = \emptyset$. Suppose by contradiction that $I \neq \emptyset$. Pick some $i \in I$ and $j \in J$ such that $a_{ij} \neq 0$ (irreducibility of $A$), and take any $k \in \mathcal{S}_i$. Lemma 3.6 gives that either $x_k > 0$ or $f_k(x) > 0$. Since $x$ is an equilibrium, $f_k(x) = 0$. So $x_k > 0$ for all $k \in \mathcal{S}_i$, contradicting the fact that $i \in I$.

[2 $\Rightarrow$ 3] Suppose that the second property holds. Pick any $j \in \{1, \ldots, m\}$. As there is some $k \in \mathcal{S}_j$ such that $x_k = 0$, (20) gives $\theta(x_1)^{b_{1j}} \theta(x_2)^{b_{2j}} \ldots \theta(x_n)^{b_{nj}} = 0$.

[3 $\Rightarrow$ 1] Obvious. ∎

## 3.5 Invariance Proofs

The key technical fact is as follows.

**Lemma 3.8** Suppose that $x : [0, t^*] \to \mathbb{R}^n$ is any solution of (1) such that $x(t) \in \mathbb{R}^n_{\geq 0}$ for all $t \in [0, t^*]$. Then the following implication holds for any $k = 1, \ldots, n$:

$$x_k(0) > 0 \;\; \Rightarrow \;\; x_k(t^*) > 0\,.$$

*Proof.* Suppose that $k$ is so that $x_k(0) > 0$. We consider, for $t \in [0, t^*]$, the scalar function $y(t) := x_k(t)$ and the functions $\alpha(t) := \alpha_k(x(t))$ and $\beta(t) := \beta_k(x(t))$. By (10),

$$\dot{y}(t) \;=\; \alpha(t) \theta(y(t)) + \beta(t)\,.$$

Since $\theta$ is locally Lipschitz, $\beta(t) \geq 0$ for all $t$, and $y(0) > 0$, it follows by a routine argument on differential equations that $y(t^*) > 0$, as desired. (The argument is: Let $z$ solve $\dot{z} = \alpha(t)\theta(z)$, $z(0) = y(0) > 0$. Since $\theta$ is locally Lipschitz, and $0$ is an equilibrium of this equation, $z(t) > 0$ for all $t$ in its domain of definition. Moreover, we have that $\dot{z} = g_1(t, z)$ and $\dot{y} = g_2(t, y)$ with $g_1(t, p) \leq g_2(t, p)$ for all $p$ (because $\beta \geq 0$). By a standard comparison theorem, we know that $z(t) \leq y(t)$ for all $t$ in the common domain of definition of $z$ and $y$. Since $y(t^*)$ is well-defined, $z(t)$ remains bounded, and thus is defined as well for $t = t^*$. So, $y(t^*) \geq z(t^*) > 0$. ∎

**Remark 3.9** Note that the only facts used in the proof are that $\theta$ is locally Lipschitz and that $\theta(y) \geq 0$ for $y \geq 0$ (which implies $\beta \geq 0$); $\theta(y) \geq 0$ for all $y$ is never needed. We assumed $\theta \geq 0$ just for convenience when displaying the general form of the systems being considered. □

**Corollary 3.10** The set $\mathbb{R}^n_+$ is forward invariant for (1).

*Proof.* Consider any solution $x : [0, T] \to \mathbb{R}^n$ of (1), and suppose that $x(0) \in \mathbb{R}^n_+$. We must prove that $x(T) \in \mathbb{R}^n_+$. Since $x(0)$ is in the interior of $\mathbb{R}^n_{\geq 0}$, the only way that the conclusion could fail is if $x(t) \in \mathbb{R}^n_0$ for some $t \in (0, T]$. We assume that this happens and derive a contradiction. Let $t^* := \min\{t \in [0, T] \,|\, x(t) \in \mathbb{R}^n_0\} > 0$. By minimality, for all $i$, $x_i(t) > 0$ for all $t \in [0, t^*)$, and in particular $x_i(t) \in \mathbb{R}^n_{\geq 0}$ for all $t \in [0, t^*]$, and also there is some index $k$ such that $x_k(t^*) = 0$. But this contradicts the conclusion of Lemma 3.8. ∎

The closure of an invariant set is also invariant, so:

**Corollary 3.11** The set $\mathbb{R}^n_{\geq 0}$ is forward invariant for (1).

**Lemma 3.12** Consider any solution $x : [0, T] \to \mathbb{R}^n$ of (1) for which $x(0) \in \mathbb{R}^n_{\geq 0}$. Suppose that there is some $j \in \{1, \ldots, m\}$ such that $x_\ell(0) > 0$ for all $\ell \in \mathcal{S}_j$. Then, $x(t) \in \mathbb{R}^n_+$ for all $t \in (0, T]$.

*Proof.* We know that $x(t) \in \mathbb{R}^n_{\geq 0}$, by Corollary 3.11. We need to see that every $i$ belongs to the set

$$I := \{i \in \{1, \ldots, m\} \,|\, x_k(t) > 0 \;\forall\, k \in \mathcal{S}_i, \forall\, t \in (0, T]\}$$

since then $\bigcup_i \mathcal{S}_i = \{1, \ldots, n\}$ gives the desired conclusion. Pick any $i$ such that $a_{ij} \neq 0$ and any $k \in \mathcal{S}_i$. Then, Lemma 3.6 says that $x_k(0) > 0$ or $f_k(x(0)) > 0$. If $x_k(0) > 0$, Lemma 3.8, applied on any subinterval $[0, t^*]$, says that $x_k(t) > 0$ for all $t$. If, instead, $f_k(x(0)) > 0$, then $\dot{x}_k(0) > 0$ and $x_k(0) \geq 0$ imply that $x_k(t) > 0$ for all $t$ small enough, so also (again by Lemma 3.8) for all $t$. Thus $i \in I$, which is therefore nonempty.

Suppose that $H := \{1, \ldots, m\} \setminus I \neq \emptyset$. Let $i \in I$ and $h \in H$ be so that $a_{hi} \neq 0$ (irreducibility of $A$). We will show that, for any given $t_0 \in (0, T]$, and for any given $k \in \mathcal{S}_h$, $x_k(t_0) > 0$, and this will contradict $h \in H$.

Since $i \in I$, $x_\ell(t_0/2) > 0$ for all $\ell \in \mathcal{S}_i$. Then, we can apply once again Lemma 3.6, which now says that $x_k(t_0/2) > 0$ or $f_k(x(t_0/2)) > 0$. As before, $x_k(t_0/2) > 0$ implies via Lemma 3.8 that $x_k(t_0) > 0$. And if $f_k(x(t_0/2)) > 0$, then $\dot{x}_k(t_0/2) > 0$ and $x_k(t_0/2) \geq 0$ imply again that $x_k(t_0) > 0$. ∎

**Corollary 3.13** Consider any solution $x : [0, T] \to \mathbb{R}_{\geq 0}^n$ of (1) for which $x(0) \notin E_0$. Then, $x(t) \in \mathbb{R}_+^n$ for all $t \in (0, T]$.

*Proof.* By Proposition 3.7, $x(0) \notin E_0$ implies that there is some $j \in \{1, \ldots, m\}$ such that $x_\ell \neq 0$ for all $\ell \in \mathcal{S}_j$. So, Lemma 3.12 insures that $x(t) \in \mathbb{R}_+^n$ for all $t \in (0, T]$. ∎

We conclude that every trajectory starting on the boundary $\mathbb{R}_0^n$ which is not an equilibrium must immediately enter the positive orthant.

# 4  Stability Proofs

We start by establishing some useful estimates.

**Lemma 4.1** Define the following quadratic function:

$$Q(\eta_1, \ldots, \eta_m) := \sum_{i=1}^m \sum_{j=1}^m a_{ij}(\eta_i - \eta_j)^2 \,. \tag{23}$$

Then, there exists a constant $\kappa > 0$ such that

$$Q(q_1, \ldots, q_m) \geq \kappa \sum_{i=1}^m \sum_{j=1}^m (q_i - q_j)^2 \tag{24}$$

for all $(q_1, \ldots, q_m) \in \mathbb{R}^m$.

*Proof.* We first observe that

$$Q(q_1, \ldots, q_m) = 0 \implies q_i = q_m, i = 1, \ldots, m - 1 \,.$$

Indeed, obviously $Q(q_1, \ldots, q_m) = 0$ implies $q_i = q_j$ for each pair $i, j$ for which $a_{ij} \neq 0$. Now let $I$ be the set of indices $i$ such that $q_i = q_m$, and $J$ its complement; as $m \in I$, $I \neq \emptyset$. We need to see that $J = \emptyset$. Suppose that $J \neq \emptyset$. The connectedness of the incidence graph of $A$ provides an $i \in I$ and $j \in J$ such that $a_{ij} \neq 0$. Thus, $q_j = q_i = q_m$, contradicting $j \in J$.

Now consider the following quadratic form in $m - 1$ variables:

$$P(\xi_1, \ldots, \xi_{m-1}) := \sum_{i=1}^{m-1} \sum_{j=1}^{m-1} a_{ij}(\xi_i - \xi_j)^2 + \sum_{i=1}^{m-1} a_{im}\xi_i^2 + \sum_{j=1}^{m-1} a_{mj}\xi_j^2 \,.$$

Since $(\eta_i - \eta_m) - (\eta_j - \eta_m) = \eta_i - \eta_j$ for all $i, j$, one has

$$Q(\eta_1, \ldots, \eta_m) = P(\eta_1 - \eta_m, \ldots, \eta_{m-1} - \eta_m) \,.$$

Note that $P$ is positive definite: if $P(q_1, \ldots, q_{m-1}) = 0$, then $Q(q_1, \ldots, q_{m-1}, 0) = 0$, which as already observed implies that all $q_i = 0$. Thus, there is some constant $\kappa_0 > 0$ such that

$$P(p_1, \ldots, p_{m-1}) \geq \kappa_0 \sum_{i=1}^{m-1} p_i^2$$

for all $(p_1, \ldots, p_{m-1}) \in \mathbb{R}^{m-1}$, which means that

$$Q(q_1, \ldots, q_m) \;\geq\; \kappa_0 \sum_{i=1}^{m-1} (q_i - q_m)^2 \tag{25}$$

for all $(q_1, \ldots, q_m) \in \mathbb{R}^m$. As $(q_i - q_j)^2 \leq 2(q_i - q_m)^2 + 2(q_j - q_m)^2$ for all $i, j$, we may re-express the estimate (25) in the form (24), using a smaller constant $\kappa$ which depends only on $\kappa_0$ and $m$. ∎

The following estimate will be the basis of a Lyapunov function property to be established later.

**Lemma 4.2** There exist two continuous functions

$$v \,:\, \mathbb{R}^n \to \mathbb{R}^n, \;\; c \,:\, \mathbb{R}_+^n \to \mathbb{R}_+$$

such that, for every pair of points $x, z$ in $\mathbb{R}_+^n$:

$$\langle \vec{\rho}(x) - \vec{\rho}(z), f(x) \rangle \;\leq\; -c(z)\delta(x, z) \,+\, \langle v(\vec{\rho}(x) - \vec{\rho}(z)), f(z) \rangle. \tag{26}$$

*Proof.* As $B$ has full column rank, there is an $m \times n$ matrix $B^{\#}$ (for instance, its pseudo-inverse) such that $B^{\#}B = I$. We let

$$v(\sigma_1, \ldots, \sigma_n) \;:=\; \left( (e^{\langle b_j, \sigma_1 \rangle}, \ldots, e^{\langle b_j, \sigma_m \rangle}) B^{\#} \right)'$$

and

$$c_0(\zeta) \;:=\; \min_{j=1,\ldots,m} e^{\langle b_j, \vec{\rho}(\zeta) \rangle}.$$

Now take any pair of positive vectors $x, z$. Denote, for each $j = 1, \ldots, m$:

$$q_j \;:=\; \langle b_j, \vec{\rho}(x) - \vec{\rho}(z) \rangle$$

and observe that

$$\langle b_i, v(\vec{\rho}(x) - \vec{\rho}(z)) \rangle \;=\; e^{q_i}, \;\; i = 1, \ldots, m$$

so, using formula (7),

$$\langle v(\vec{\rho}(x) - \vec{\rho}(z)), f(z) \rangle \;=\; \sum_{i=1}^{m} \sum_{j=1}^{m} a_{ij} e^{\langle b_j, \vec{\rho}(z) \rangle} \left( e^{q_i} - e^{q_j} \right). \tag{27}$$

Therefore (writing $g(x, z) = \langle v(\vec{\rho}(x) - \vec{\rho}(z)), f(z) \rangle$ for simplicity):

$$
\begin{aligned}
\langle \vec{\rho}(x) - \vec{\rho}(z), f(x) \rangle 
&= \sum_{i=1}^{m} \sum_{j=1}^{m} a_{ij} e^{\langle b_j, \vec{\rho}(x) \rangle} (q_i - q_j) \\
&= \sum_{i=1}^{m} \sum_{j=1}^{m} a_{ij} e^{\langle b_j, \vec{\rho}(z) \rangle} e^{q_j} (q_i - q_j) \\
&= \sum_{i=1}^{m} \sum_{j=1}^{m} a_{ij} e^{\langle b_j, \vec{\rho}(z) \rangle} (e^{q_j}(q_i - q_j) - (e^{q_i} - e^{q_j})) \,+\, g(x, z) \\
&\leq -\frac{1}{2} \sum_{i=1}^{m} \sum_{j=1}^{m} a_{ij} e^{\langle b_j, \vec{\rho}(z) \rangle} (q_i - q_j)^2 \,+\, g(x, z) \\
&\leq -\frac{1}{2} c_0(z) \sum_{i=1}^{m} \sum_{j=1}^{m} a_{ij} (q_i - q_j)^2 \,+\, g(x, z) \\
&= -\frac{1}{2} c_0(z)\, Q(q_1, \ldots, q_m) \,+\, g(x, z),
\end{aligned}
$$

$$\tag{28}$$

$$\tag{29}$$

where $Q$ is the quadratic form in Lemma 4.1. Equality (28) follows by adding and subtracting $g(x, z)$ and using (27). To justify (29), we note first that, for each $a > 0$, the function $\mathbb{R}_{\geq 0} \to \mathbb{R}$:

$$f_a(r) := e^a(r-a) - e^r + e^a + \frac{1}{2}(r-a)^2$$

is always $\leq 0$ (because $f_a(0) = -e^a a - 1 + e^a + (1/2)a^2 < 0$, $f_a(r) \to -\infty$ as $r \to +\infty$, and $f_a'(r) = e^a - e^r + (r-a) \neq 0$ for all $r > 0$). Now we use the inequality $e^a(r-a) - e^r + e^a \leq -\frac{1}{2}(r-a)^2$ in each term of the sum with $a = q_j$ and $r = q_i$ (recall $a_{ij} e^{\langle b_j, \vec{\rho}(z) \rangle} \geq 0$).

Lemma 4.1 gives that $Q(q_1, \ldots, q_m) \geq \kappa\delta(x, z)$. Thus, we may take $c(z) := \kappa c_0(z)/2$. ∎

## 4.1  An Entropy Distance

Recall that we are assuming that $\int_0^1 |\rho(r)|\, dr < \infty$. For any fixed constant $c \in \mathbb{R}$, we consider the following function:

$$R_c(r) := \int_1^r \rho(s)\, ds - cr\,.$$

This function is a well-defined continuous mapping $\mathbb{R}_{\geq 0} \to \mathbb{R}$, continuously differentiable for $r > 0$. Moreover, $R_c$ is strictly convex, since for $r > 0$ its derivative $\rho(r) - c$ is strictly increasing and onto $\mathbb{R}$; it achieves a global minimum at the unique $r_c \in \mathbb{R}_+$ where $\rho(r_c) = c$, decreases for $r \in [0, r_c]$, and increases to $+\infty$ for $r > r_c$.

The following function will play a central role:

$$W : \mathbb{R}_{\geq 0}^n \times \mathbb{R}_+^n : (x, z) \mapsto \sum_{i=1}^n R_{\rho(z_i)}(x_i)\,.$$

The above-mentioned properties of the functions $R_{\rho(z_i)}$ imply that

$$x \neq z \quad \Rightarrow \quad W(x, z) > W(z, z)\,, \tag{30}$$

i.e., for each fixed $z \in \mathbb{R}_+^n$, the function $W(\cdot, z)$ has a unique global minimum, at $z$. Note also that the gradient of $W(\cdot, z)$:

$$\frac{\partial W}{\partial x}(x, z) = (\vec{\rho}(x) - \vec{\rho}(z))' \tag{31}$$

(defined for $x \in \mathbb{R}_+^n$) vanishes only at $x = z$ and that (since $R_{\rho(z_i)}(x_i) \to +\infty$ if $x_i \to +\infty$)

$$|x| \to +\infty \quad \Rightarrow \quad W(x, z) \to +\infty \tag{32}$$

for every given $z$. As $W(\cdot, z)$ is continuous, this implies that

$$\{x \mid W(x, z) \leq w\} \tag{33}$$

is compact for every $z$ and every $w \in \mathbb{R}$.

**Remark 4.3** In the special case $\rho = \ln$, $W(x, z) = \sum_{i=1}^n x_i \ln x_i - x_i - x_i \ln z_i$. Then this formula, when states $x$ are interpreted probabilistically in applications such as chemical networks, is suggested by "relative entropy" considerations. □

## 4.2  Main Stability Results

For any $\xi \in \mathbb{R}_+^n$, let us denote by $S^\xi$ the positive class $(\xi + \mathcal{D}) \bigcap \mathbb{R}_{\geq 0}^n$ containing $\xi$, $\bar{x}^\xi$ the unique interior equilibrium $\bar{x} \in S^\xi$, $E_0^\xi$ the set (possibly empty) of boundary equilibria in $S^\xi$, i.e., the set $S^\xi \bigcap E_0$, and

$$E^\xi = \{\bar{x}^\xi\} \bigcup E_0^\xi = S^\xi \bigcap E$$

the set of all equilibria in $S^\xi$. The main technical fact, which will imply the completeness Lemma 2.2 as well as the convergence parts in Theorems 1 and 2, is as follows.

**Proposition 4.4** For each $\xi \in \mathbb{R}_{\geq 0}^n$, the solution $x(\cdot)$ of the initial value problem $\dot{x} = f(x)$ with $x(0) = \xi$ is defined for all $t \geq 0$, and

$$x(t) \to E^\xi \ \text{ as } \ t \to +\infty \,.$$

*Proof.* Fix $\xi$ in $\mathbb{R}_{\geq 0}^n$. We will use

$$V(x) \ := \ W(x, \bar{x}^\xi) - W(\bar{x}^\xi, \bar{x}^\xi)$$

as a Lyapunov-like function. By (30), this function is positive definite relative to the equilibrium $\bar{x}^\xi$, i.e., $V(x) \geq 0$ for all $x \in \mathbb{R}_{\geq 0}^n$, and $V(x) = 0$ if and only if $x = \bar{x}^\xi$. Moreover, $V$ is proper, meaning that the sublevel sets $\{x \mid V(x) \leq w\}$ are compact, for all $w \in \mathbb{R}_{\geq 0}$, by (33). Finally, $V$ is continuously differentiable in the interior $\mathbb{R}_+^n$, and

$$\nabla V(x) f(x) \ \leq \ -c\,\delta(x, \bar{x}^\xi) \tag{34}$$

for all $x \in \mathbb{R}_+^n$, where $c = c(\bar{x}^\xi) > 0$, by (31) and (26). In particular, since $\delta(x, \bar{x}^\xi) = 0$ implies that $x$ is an equilibrium, and since there is a unique interior equilibrium in each class:

$$x \in S^\xi \textstyle\bigcap \mathbb{R}_+^n, \ x \neq \bar{x}^\xi \ \Rightarrow \ \nabla V(x) f(x) < 0 \,. \tag{35}$$

Now consider the maximal solution $x(\cdot)$ starting from $x(0) = \xi$, which is a priori defined on some interval $[0, t^*)$. If $\xi \in E$, there is nothing to prove, so we assume that $\xi$ is not an equilibrium. As, in particular, $\xi \notin E_0$, Corollary 3.13 insures that $x(t) \in \mathbb{R}_+^n$ for all $t > 0$ in the domain of definition. So $V(x(t))$ is differentiable for $t > 0$, and $dV(x(t))/dt \leq 0$ (by (34)), which means that $V(x(t))$ is nondecreasing. Since $V$ is proper, this means that the maximal trajectory is precompact, and hence it is defined on the entire interval $[0, +\infty)$. Furthermore, the LaSalle Invariance Theorem implies that

$$x(t) \to M \ \text{ as } \ t \to +\infty \,,$$

where $M$ is the largest invariant set included in $\{p \mid V(p) = a\}$, for some $a \geq 0$.

Since $x(t) \in S^\xi$ for all $t$, and $S^\xi$ is closed, any limit point of $x(\cdot)$ is included in $S^\xi$ as well. So $M \subseteq S^\xi$. Thus, we need only to show that $M \subseteq E^\xi$. Pick any $\zeta \in M$, and consider the forward trajectory $z(t)$ starting from $\zeta$.

Suppose first that $\zeta \in \mathbb{R}_+^n$. Assume that $\zeta \neq \bar{x}^\xi$. Then (35) says that $dV(z(t))/dt < 0$ at $t = 0$, which implies that $V(z(t)) < a$ for $t$ small, contradicting the fact that $M$ is included in $\{p \mid V(p) = a\}$. Thus $\zeta = \bar{x}^\xi$.

Suppose now $\zeta \in \mathbb{R}_0^n$. If $\zeta \notin E_0$, then Corollary 3.13 gives that $z(t_0) \in \mathbb{R}_+^n$ for some $t_0 > 0$. But invariance of $M$ says that $z(t_0) \in M$, and we already showed that $M \bigcap \mathbb{R}_+^n \subseteq E$. So, $z(t_0)$ is an equilibrium, which means that $z(t) \equiv z(t_0)$ and hence $\zeta = z(t_0) \notin \mathbb{R}_0^n$, a contradiction again. ∎

**Remark 4.5** Observe that the proof actually shows that the domain of attraction of $\bar{x}^\xi$ contains $\{x \mid V(x) < w_0\}$, where $w_0 := \min_{x \in E_0 \bigcap S^\xi} V(x)$. □

Theorem 1 follows from the above. To prove Theorem 2, part *b*, pick any positive class $S$ and $\xi \in S$. In the above proof, we simply note that (35) holds for all $x$ in some neighborhood of $\bar{x}^\xi = \bar{x}_S$, so $V$ is a (local) Lyapunov function for the system restricted to the class $S$. Finally, to see Theorem 2, part *c*, note that if $\zeta \in E_0 \bigcap S^\xi \neq \emptyset$ then $\zeta$ (being an equilibrium) is not attracted to $\bar{x}^\xi$, and if instead $E_0 \bigcap S^\xi = \emptyset$ then the Proposition tells us that all trajectories converge to $E^\xi = \{\bar{x}^\xi\}$.

# 5  The Kinetic Proofreading Example, and Chemical Networks

We now show the global stability and unique equilibrium properties of the kinetic proofreading model which was described in the introduction. For this purpose, we wish to see the equations as those of an appropriate system (1), when restricted to a suitable class (which is determined by the constants $M^*$ and $T^*$). Recalling the conservation laws $M + \sum C_i = M^*$ and $T + \sum C_i = T^*$, we write the equations as a system of dimension $n = N + 3$:

$$
\begin{aligned}
\dot{T} &= -k_1 T M + \sum_{i=0}^{N} k_{-1,i} C_i \\
\dot{M} &= -k_1 T M + \sum_{i=0}^{N} k_{-1,i} C_i \\
\dot{C}_0 &= k_1 T M - (k_{-1,0} + k_{p,0}) C_0 \\
&\vdots \\
\dot{C}_i &= k_{p,i-1} C_{i-1} - (k_{-1,i} + k_{p,i}) C_i \\
&\vdots \\
\dot{C}_N &= k_{p,N-1} C_{N-1} - k_{-1,N} C_N .
\end{aligned}
$$

This is indeed a system of form (1). To see this, we use $x = (T, M, C_0, \ldots, C_N)'$ as a state, and take $m = n - 1 = N + 2$,

$$
b_1 = \begin{pmatrix} 1 \\ 1 \\ 0 \\ 0 \\ \vdots \\ 0 \end{pmatrix}, \; b_2 = \begin{pmatrix} 0 \\ 0 \\ 1 \\ 0 \\ \vdots \\ 0 \end{pmatrix}, \; b_3 = \begin{pmatrix} 0 \\ 0 \\ 0 \\ 1 \\ \vdots \\ 0 \end{pmatrix}, \; \ldots \; b_m = \begin{pmatrix} 0 \\ 0 \\ 0 \\ 0 \\ \vdots \\ 1 \end{pmatrix},
$$

and $A = (a_{ij})$ with $a_{21} = k_1$, $a_{1i} = k_{-1,i-2}$ $(i = 2, \ldots, m)$, $a_{i,i-1} = k_{p,i-3}$ $(i = 3, \ldots, m)$, and all other $a_{ij} = 0$.

Thus, $\mathcal{D} = \{x \mid T + C_0 + \ldots + C_N = M + C_0 + \ldots + C_N = 0\}$, and the positive classes are of the form $S = S^{\alpha,\beta}$, intersections with $\mathbb{R}_{\geq 0}^n$ of the affine planes

$$
T + C_0 + \ldots + C_N = \alpha , \quad M + C_0 + \ldots + C_N = \beta ,
$$

with $\alpha > 0$ and $\beta > 0$. The original system is nothing else than the class determined by $\alpha = T^*$ and $\beta = M^*$. Thus, the conclusions will follow from Theorem 2 as soon as we prove that $S \bigcap E_0 = \emptyset$ for any positive class $S$. To see this, we may use Proposition 3.7. Pick any $x \in \mathbb{R}_0^n$ and any positive class $S^{\alpha,\beta}$. We must find some $j \in \{1, \ldots, m\}$ with the property that $x_k \neq 0$ for all $k \in \mathcal{S}_j$. In our application, $\mathcal{S}_1 = \{1, 2\}$ and $\mathcal{S}_j = \{j + 1\}$ for $j = 2, \ldots, m$. If the property is not satisfied for some $j \in \{2, \ldots, m\}$, then $C_i = 0$ for all $i$. But in this case, the equations for $S^{\alpha,\beta}$ give that $T = \alpha > 0$ and also $M = \beta > 0$, so neither can $j = 1$ be used. In conclusion, $x \notin E_0$, and hence part $c$ of the theorem applies. ∎

## 5.1 Comments on Chemical Networks

Let us very briefly indicate the motivation for the form of the systems (1) as arising from mass-action kinetics in chemical network theory. One studies $n$ chemical species $A_1, \ldots, A_n$, whose concentrations as a function of time are given by $x_1(t), \ldots, x_n(t)$ and then derives a system of differential equations for the $x_i$'s on the basis of the known reactions that occur among the substances $A_i$. For instance, suppose that each molecule of $A_1$ can react with four molecules of $A_2$ to produce two molecules of $A_3$. This is indicated graphically by

$$
A_1 + 4A_2 \; \to \; 2A_3 .
$$

Assuming that the reactor is well-mixed, the probability of such a reaction occurring, at an instant $t$, is proportional to the product $x_1(t) x_2(t)^4$. Since we gain two molecules of $A_3$ for each such reaction, this gives rise to a rate of increase

$$
\dot{x}_3 = 2k x_1 x_2^4 \tag{36}
$$

for the concentration of $A_3$, where $k$ is a suitable constant of proportionality. This constant $k$ is often thought of as a "reaction rate" and one writes graphically:

$$A_1 + 4A_2 \xrightarrow{k} 2A_3 \,.$$

In this manner, one puts together the whole system of differential equations. A convenient way to specify the resulting system is by building a matrix $A$ from the reaction rates, and introducing vectors $b_1, \ldots, b_m$ to describe each of the "complexes" such as $A_1 + 4A_2$ and $2A_3$ by specifying the contributions from each type of molecule. For example, $A_1 + 4A_2$ might give rise to the vector $b_1 = (1, 4, 0, 0, \ldots, 0)'$ and $2A_3$ to the vector $b_2 = (0, 0, 2, 0, \ldots, 0)'$, and reaction (36) contributes then the term $kx_1x_2^4(b_2 - b_1)$ to the differential equations (1), where the reason for "$kx_1x_2^4$" is obvious, and $b_2 - b_1$ corresponds to the fact that the reaction goes from the complex associated to the vector $b_1$ to that one associated to $b_2$. (Note how the equation for $\dot{x}_3$ will then have the contribution $kx_1x_2^4(b_{32} - b_{31}) = 2kx_1x_2^4$, and also $\dot{x}_1$ will have a term $kx_1x_2^4(b_{12} - b_{11}) = -kx_1x_2^4$ to indicate that $A_1$ is being eliminated at the given rate.) In this context, the space $\mathcal{D}$ is called the *stoichiometric subspace* associated to the reaction, and a class is a *stoichiometric compatibility class*.

The results explained in Theorems 1 and 2 are basically the main theorems for what are called mass-action networks of zero deficiency and a single linkage class. We do not define these terms here. For more details, see for instance [1, 2]. (A small remark for readers who compare our results with those in the chemical network literature: Condition (3) might appear to be slightly stronger than needed, since the zero-deficiency theorem would only require the columns $b_i$ to be affinely, rather than linearly, independent. However, no generality is lost, because if we start with an affinely independent set of vectors $v_i$, we can introduce the vectors $b_i = (1, v_i')'$ ( i.e, just add a constant coordinate $= 1$) in one more dimension. The span of the differences $v_i - v_j$ has the same dimension as $\mathcal{D}$. So we need only consider a new set of differential equations in which we add a variable satisfying $\dot{x}_0 \equiv 0$, to bring everything into the setup considered in this note.)

# References

[1] Feinberg, M., "Chemical reaction network structure and the stabiliy of complex isothermal reactors - I. The deficiency zero and deficiency one theorems," Review Article 25, *Chemical Engr. Sci.* **42**(1987): 2229-2268.

[2] Horn, F.J.M., and Jackson, R., "General mass action kinetics," *Arch. Rational Mech. Anal.* **49**(1972): 81-116.

[3] McKeithan, T.W., "Kinetic proofreading in T-cell receptor signal transduction," *Proc. Natl. Acad. Sci. USA* **92**(1995): 5042-5046.

[4] Wofsy, C., "Modeling receptor aggregation," series of lectures at the *Workshop on Mathematical Cellular Biology*, Pacific Institute for the Mathematical Sciences, UBC, Vancouver, August 1999.

[5] Sontag, E.D., in preparation.